\documentclass[a4paper,12pt]{amsart}

\usepackage{amsmath,amssymb,amsthm,latexsym}

\newtheorem*{theorem}{Theorem}
\newtheorem*{corollary}{Corollary}
\newtheorem*{definition}{Definition}

\newcommand{\sa}{_{\textit{\normalsize{sa}}}}

\begin{document}

\title[A Monotone Selection Principle]{A Monotone Selection Principle in C*-algebras}

\author[M. Thill]{Marco Thill}

\address{bd G.\,-\,D.\ Charlotte 53 \\
                 L\ -\,1331 Luxembourg\,-\,City}

\email{math@pt.lu}

\date{}



\begin{abstract}
We apply to operator algebra theory a monotone selection principle which apparently
escaped attention (of operator algebra theorists) so far. This principle relates to the basic
order theoretic characterisation of von Neumann algebras given by Kadison, and the
simplified form this result takes in separable Hilbert spaces. In the separable case we need
only consider increasing sequences rather than increasing nets. We apply an argument of
Klaus Floret to show that, within the realm of commutativity, there exists a general monotone
selection principle providing for this simplification. Thereby we obtain a valuable shortcut
and a handy tool for related purposes. Actually, a more general selection principle is
proved within the framework of vector lattices.
\end{abstract}

\maketitle


\section{Introduction}

We first would like to put into perspective the needed elements of operator algebra
theory. The following three notions are cornerstones of the order theoretic aspect of
operator algebras: \\
- \emph{monotone complete} hermitian parts of C*-algebras, \\
- \emph{monotone closed} subsets of the above, \\
- \emph{normal} positive linear maps between the above. \\
We shall review these notions and their relevance in the next paragraph.

Let $A$ be a C*-algebra. We shall denote by $A\sa$ the hermitian part of $A$. One
says that $A\sa$ is \emph{monotone complete} if each non-empty upper bounded
upward directed subset of $A\sa$ has a supremum in $A\sa$, cf.\ \cite[3.9.2]{Ped}.
For example the hermitian part of a von Neumann algebra is monotone complete.
Similarly, if $B$ is a C*-subalgebra of a C*-algebra $A$ with monotone complete
hermitian part $A\sa$, then $B\sa$ is called \emph{monotone closed} in $A\sa$, if
$B\sa$ contains the supremum in $A\sa$ of each non-empty upward directed subset
of $B\sa$ which is upper bounded in $A\sa$. The following theorem of Kadison
occupies a central place in the order theoretic part of operator algebra theory. If
$H$ is a Hilbert space, and $A$ is a C*-subalgebra of $B(H)$, then $A$ is a von
Neumann algebra if and only if $A\sa$ is monotone closed in $B(H)\sa$. For a proof,
see \cite[2.4.4]{Ped}. Let $\phi : A \to B$ be a positive linear map between two
C*-algebras $A$ and $B$. (For example $\phi$ a state or a representation.) If $A\sa$
and $B\sa$ are monotone complete, then $\phi$ is called \emph{normal}, if for each
non-empty upper bounded upward directed subset $J$ of $A\sa$ one has
$\phi (\sup J) = \sup _{j \in J} \phi (j)$, cf.\ \cite[2.5.1]{Ped}. It can be shown that if
$\phi : A \to B$ is a normal $*$-algebra homomorphism between C*-algebras $A$
and $B$ with monotone complete hermitian parts, then $\phi (A)\sa$ is a monotone
closed C*-subalgebra of $B\sa$, cf.\ the proof of \cite[2.5.3]{Ped}. Thus, if $B$
moreover is a von Neumann algebra, then $\phi (A)$ is a von Neumann algebra.

It has been noted that on a separable Hilbert space, the above conditions can be
relaxed. Namely, it is enough to consider increasing sequences instead of increasing
nets, cf.\ \cite[2.4.3]{Ped}. One then speaks of: \\
- \emph{monotone sequentially complete} hermitian parts of C*algebras, \\
- \emph{monotone sequentially closed} subsets of the above, \\
- \emph{sequentially normal} positive linear maps between the above. \\
It has also been noted that for many purposes, the assumption of separability of the
Hilbert space can be relaxed to assuming that the von Neumann algebra is
$\sigma$-finite. We shall review this notion in the next three paragraphs.

A state $\psi$ on a C*-algebra $A$ is called \emph{faithful} if $\psi (a) > 0$ for every
$a \in A _+ \setminus \{ 0 \}$, cf.\ \cite[3.7.2]{Ped}. The following fact is well-known:
in a C*-algebra carrying a faithful state every set of pairwise orthogonal projections is
countable, cf.\ \cite[2.5.6]{BR}. It is to this well-known fact that we shall adjoin a versatile
supplement. 

(We open a parenthesis for a proof here because it is so short. Let $\psi$ be a faithful state
on a C*-algebra $A$. Let $(\pi _\psi, H _\psi, c _\psi)$ be the GNS representation of $\psi$.
Let $P$ be a set of pairwise orthogonal projections in $A$. Then $(\pi _\psi (p) ) _{p \in P}$
is a family of pairwise orthogonal projections in $H _\psi$. Let $q$ denote its sum. By the
Parseval equality it follows that
$\| q c _\psi \| ^2 = \sum _{p \in P} \| \pi _\psi (p) c _\psi \| ^2$.
Since the left side of this equality is finite, the sum on the right can only contain countably
many non-zero terms. All other terms satisfy $0 = \| \pi _\psi (p) c _\psi \| ^2 = \psi (p^*p)$,
whence $0 = p^*p = p$ as $\psi$ is faithful.)

A von Neumann algebra $\mathfrak{M}$ is called \emph{$\sigma$-finite} if every set of
pairwise orthogonal projections in $\mathfrak{M}$ is countable, cf.\ \cite[2.5.1]{BR},
or \cite[3.8.3]{Ped}. Obviously a Hilbert space $H$ is separable if and only if the von Neumann
algebra of all bounded operators on $H$ is $\sigma$-finite. Also, a von Neumann algebra is
$\sigma$-finite if and only if it admits a faithful state. This follows from the above proof together
with \cite[2.5.6]{BR}.  Note that $\sigma$-finiteness of a von Neumann algebra $\mathfrak{M}$
is of a \emph{commutative} nature in the sense that it pertains to commutative C*-subalgebras
of $\mathfrak{M}$ only.

We give a related selection principle which apparently escaped attention so far. The
argument is due to Klaus Floret, cf.\ his proof \cite[p.\ 228]{Fl} of the order completeness
of the $L^p$ spaces for $1 \leq p < \infty$. A core statement goes as follows. (For a more
precise version with proof, see the next section.)

\begin{theorem}\label{core}
Let $A$ be a commutative C*-algebra carrying a faithful state $\psi$. If $A\sa$ is
monotone sequentially complete, then it is is monotone complete, and the following
\emph{monotone selection principle} holds. For every non-empty upper bounded
upward directed subset $J$ of $A\sa$ there exists an increasing sequence $(j_n)$
in $J$ with the same supremum as $J$. Every increasing sequence $(j_n)$ in $J$
with $\sup _n \psi (j_n) = \sup _{j \in J} \psi (j)$ does the job.
\end{theorem}

We feel that this principle is rather smooth in operation compared to the clumsy notion
of $\sigma$-finiteness. We have the following consequences.

\begin{corollary}
Let $\phi : A \to B$ be a positive linear map between two C*-algebras $A$ and $B$
with monotone complete hermitian parts. Assume that $A$ is commutative and admits
some faithful state. If $\phi$ is sequentially normal, then $\phi$ is normal.
\end{corollary}

\begin{corollary}
Let $A$ be a commutative C*-algebra carrying a faithful state $\psi$. If $A\sa$ is
monotone sequentially complete, and if $\psi$ is sequentially normal, then $A\sa$
is monotone complete, and $\psi$ is normal. In this case, the GNS representation
$\pi _\psi$ is a normal isomorphism of C*-algebras from $A$ onto a von Neumann
algebra.
\end{corollary}

\begin{proof}
The proof of the statement concerning the GNS representation follows from
the proofs of \cite[3.3.9]{Ped} and \cite[2.5.3]{Ped}.
\end{proof}

In the next section, a more general selection principle is formulated and proved in the abstract
framework of vector lattices. This makes that commutative C*-algebras belong within the reach
of the result. Indeed, it is well-known that the hermitian part of a commutative C*-algebra is a
vector lattice. (This follows for example from the commutative Gelfand-Na\u{\i}mark theorem.)
Conversely, the hermitian part of a non-commutative C*-algebra is never a vector lattice. (This
follows from the remark \cite[1.4.9]{Ped}.) However, our results are of immediate interest to
non-commutative C*-algebras as well, namely via maximal commutative C*-subalgebras.
We mention in this respect the notion of \emph{AW*-algebras} and the notion of
\emph{completely additive} positive linear maps between AW*-algebras, cf.\ \cite[3.9.2]{Ped},
and \cite[3.9.6]{Ped}. Note again that $\sigma$-finiteness is of a commutative nature as well.

\section{A more precise selection principle and proof}

In the following definition we formalise and clarify completeness and closedness with respect
to monotone nets and sequences. Also, the notion of a \emph{faithful} positive linear functional
on an ordered vector space is introduced.

\begin{definition}
Let $W$ be an ordered vector space, and let $V$ be a vector subspace of $W$.
(Typically $W$ the hermitian part of a von Neumann algebra.)

We shall say that $W$ is \emph{monotone complete}, if each non-empty upper
bounded upward directed subset of $W$ has a supremum in $W$. We shall say that $W$
is \emph{monotone sequentially complete}, if each upper bounded increasing sequence
in $W$ has a supremum in $W$.

If $W$ is monotone complete, then $V$ shall be called \emph{monotone closed}
in $W$, if $V$ contains the supremum $\sup J$ in $W$ of each non-empty upward directed
subset $J$ of $V$ that is upper bounded in $W$. (It is clear that then $\sup J$ also is
the supremum in $V$ of $J$, and that $V$ then is monotone complete.) If $W$ is monotone
sequentially complete, then $V$ shall be called \emph{monotone sequentially closed} in $W$,
if $V$ contains the supremum $\sup _n j_n$ in $W$ of each increasing sequence $(j_n)$
in $V$ that is upper bounded in $W$. (It is clear that then $\sup _n j_n$ also is the
supremum in $V$ of $(j_n)$, and that $V$ then is monotone sequentially complete.)

If $W _+$ denotes the set of positive elements of $W$, then a positive linear functional
$\psi$ on $W$ shall be called \emph{faithful}, if $\psi (a) > 0$ for each
$a \in W _+ \setminus \{ 0 \}.$
\end{definition}

\begin{theorem}
Let $W$ be a vector lattice carrying a faithful positive linear functional $\psi$. Assume
that $W$ is monotone sequentially complete, and that $V$ is a monotone sequentially
closed vector subspace of $W$. Then $W$ is monotone complete, and $V$ is monotone
closed in $W$.

Indeed, which is more, the following \emph{monotone selection principle} holds.
Whenever $J$ is a non-empty upward directed subset of $V$ that is upper bounded in $W$,
there exists an increasing sequence $(j_n)$ in $J$ with
\[ \sup _n j_n = \sup J \text{ in } W. \]
Every increasing sequence $(j_n)$ in $J$ with
\[ \sup _n \psi (j_n) = \sup _{j \in J} \psi (j) \]
does the job. From the result that $V$ is monotone closed in $W$, one also has
\[ \sup _n j_n = \sup J \text{ in } V. \]
\end{theorem}

\begin{proof} Let $J$ be a non-empty upward directed subset of $V$ that is upper bounded in
$W$. Let $(j_n)_{n \geq 1}$ be any increasing sequence in $J$ with
\[ \sup _{n \geq 1} \psi (j_n) = \sup _{j \in J} \psi (j). \tag*{$(*)$} \]
We define
\[ j_0 = \sup _{n \geq 1} j_n \text{ in } W, \]
which is possible by the assumption that $W$ is monotone sequentially complete. One
has that $j_0$ belongs to $V$ by the assumption that $V$ is monotone sequentially
closed in $W$. To prove the theorem, it is enough to show that $j_0 = \sup J$ in $W$.
(This is so because if we choose $V := W$, we obtain that $W$ is monotone complete.
We clearly then also have that $V$ is monotone closed in $W$.) To prove that
$j_0 = \sup J$ in $W$, it is sufficient to show that $j_0$ is an upper bound of $J$.
(Indeed, if $j_0$ is an upper bound of $J$, and $h$ is any other upper bound of $J$
in $W$, then $h$ also is an upper bound of $(j_n)_{n \geq 1}$. Since $j_0$ is the least
upper bound of $(j_n)_{n \geq 1}$ in $W$, it follows that $h \geq j_0$. This says that
$j_0=\sup J$ in $W$.) Let $j \in J$ be arbitrary. We have to prove that $j_0 \geq j$.
For arbitrary $n \geq 1$, there exists an element $g_n$ of $J$ with
$g_n \geq j_n, g_n \geq j$ by the upward direction of $J$. One then has
\[  j - j_0 \leq g_n - j_0 \leq g_n-j_n. \]
Using that $W$ is a vector lattice, we have
\[ (j - j_0)_+ = 0 \vee (j - j_0) \leq 0 \vee (g_n - j_n) = (g_n - j_n)_+ = g_n - j_n, \]
and so
\[ \psi ((j - j_0)_+) \leq \psi (g_n - j_n)= \psi (g_n) - \psi (j_n) \leq \sup _{j \in J} \psi (j) - \psi (j_n). \]
Since $n \geq 1$ is arbitrary, it follows with $(*)$ that
\[ \psi ((j - j_0)_+) \leq \sup _{j \in J} \psi (j) - \sup _n \psi (j_n) = 0, \]
whence $(j  -j_0)_+=0$ as $\psi$ is faithful. This implies $j - j_0 \leq 0$, or $j_0 \geq j$.
\end{proof}

\begin{definition}
Let $\phi : W_1 \to W_2$ be a positive linear map between ordered vector spaces
$W_1$, $W_2$. If $W_1$, $W_2$ are monotone complete, then $\phi$ is called
\emph{normal}, if for each non-empty upper bounded upward directed subset $J$ of $W_1$,
one has $\phi (\sup J) = \sup _{j \in J} \phi (j)$. If $W_1$, $W_2$ are monotone sequentially
complete, then $\phi$ is called \emph{sequentially normal}, if for each upper bounded
increasing sequence $(j_n)$ in $W_1$, one has $\phi (\sup _n j_n) = \sup _n \phi (j_n)$.
\end{definition}

\begin{corollary}
Let $\phi : A \to B$ be a positive linear map between monotone complete ordered
vector spaces $A$ and $B$. Assume that $A$ is a vector lattice admitting some faithful
positive linear functional. If $\phi$ is sequentially normal, then $\phi$ is normal.
\end{corollary}


\begin{thebibliography}{9}
\bibitem{BR} O.\ Bratteli, D.\ Robinson, \emph{Operator Algebras and
               Quantum Statistical Mechanics I}, 2nd ed., TMP (Springer, 1987).
\bibitem{Fl}  K.\ Floret, \emph{Ma\ss- und Integrationstheorie}
               (Teubner, Stuttgart 1981).
\bibitem{Ped} G.\ K.\ Pedersen, \emph{C*-algebras and their automorphism
               groups}, LMS Monographs  No.\ 14 (Academic Press, 1979).
\end{thebibliography}
\end{document}